\documentclass[11pt]{article}

\usepackage{amsfonts,amssymb,amsmath,amsthm,epsfig,euscript}
\usepackage{epic, eepic, graphics}

\newcommand\p{\circle*{0.3}}

\newtheorem{definition}{Definition}

\newtheorem{remark}{Remark}

\newtheorem{problem}{Problem}

\long\def\symbolfootnote[#1]#2{\begingroup
\def\thefootnote{\fnsymbol{footnote}}\footnote[#1]{#2}\endgroup}

\def\mn{{\mbox{-}}}

\def\S{\mathcal{S}}

\def\PP{\mathbb{P}}
\def\E{\mathbb{E}}
\def\O{\mathbb{O}}

\newcommand{\fig}[2]{\begin{figure}[ht]
\centerline{\scalebox{.66}{\epsfig{file=#1.eps}}}
\caption{#2}
\label{figure:#1}
\end{figure}}

\newcommand{\sg}{\sigma}

\title{Place-difference-value patterns: A generalization of
generalized permutation and word patterns.}

\author{
Sergey Kitaev\footnote{The work presented here was supported by grant no. 090038011 from the Icelandic Research Fund.} \\
\small The Mathematics Institute\\[-0.8ex]
\small Reykjav\'{i}k University \\[-0.8ex]
\small IS-103 Reykjav\'{i}k, Iceland\\[-0.8ex]
\small \texttt{sergey@ru.is} \and
Jeffrey Remmel\footnote{Partially supported by NSF grant DMS 0654060.} \\
\small Department of Mathematics\\[-0.8ex]
\small University of California, San Diego\\[-0.8ex]
\small La Jolla, CA 92093-0112. USA\\[-0.8ex]
\small \texttt{remmel@math.ucsd.edu}
}

\date{\small Submitted: Date 1;  Accepted: Date 2;
 Published: Date 3.\\
\small MR Subject Classifications: 05A15}

\begin{document}
\maketitle

\begin{abstract}
Motivated by study of Mahonian statistics, in 2000, Babson and
Steingr\'{\i}msson~\cite{BS} introduced  the notion of a
``generalized permutation pattern'' (GP) which generalizes the
concept of ``classical'' permutation pattern introduced by Knuth in
1969. The invention of GPs led to a large number of publications
related to properties of these patterns in permutations and words.
Since the work of Babson and Steingr\'{\i}msson, several further
generalizations of permutation patterns have appeared in the
literature, each bringing a new set of permutation or word pattern
problems and often new connections with other combinatorial objects
and disciplines. For example, Bousquet-M\'elou et al.~\cite{BCDK}
introduced a new type of permutation pattern that allowed them to
relate permutation patterns theory to the theory of partially
ordered sets.

In this paper we introduce  yet another, more general definition of
a pattern, called place-difference-value patterns (PDVP) that covers
all of the most common definitions of permutation and/or word
patterns that have occurred in the literature. PDVPs provide many
new ways to develop the theory of patterns in permutations and
words. We shall give several examples of PDVPs in both permutations
and words that cannot be described in terms of any other pattern
conditions that have been introduced previously. Finally, we raise
several bijective questions linking our patterns to other
combinatorial objects.
\end{abstract}

\section{Introduction}

In the last decade, several hundred papers have been published on the
subject of patterns in words and permutations. This is a new, but
rapidly growing, branch of combinatorics which has its roots in the works
by Rotem, Rogers, and Knuth in the 1970s and early 1980s. However,
the first systematic study of permutation patterns was not undertaken until the paper by
Simion and Schmidt~\cite{SchSim} which appeared in 1985. The
field has experienced explosive
growth since 1992. The notion
of  patterns in permutations and words
has proven to be a useful language in a
variety of seemingly unrelated problems including the {\em theory of
Kazhdan-Lusztig polynomials}, {\em singularities of Schubert
varieties}, {\em Chebyshev polynomials}, {\em rook polynomials
for Ferrers board}, and various {\em sorting algorithms}
including {\em sorting stacks} and {\em sortable permutations}.

A (``{\em classical}'') {\em permutation pattern} is a permutation $\sg = \sg_1 \ldots \sg_k$ in the symmetric group $S_k$  viewed as a word
without repeated letters. We say that $\sg$ {\em occurs} in a
permutation $\pi = \pi_1 \ldots \pi_n$ if there is a subsequence
$1 \leq i_1 < \cdots < i_k \leq n$ such that $\pi_{i_1} \pi_{i_2} \ldots
\pi_{i_k}$ is order isomorphic to $\sg$. We say that $\pi$ {\em avoids}
$\sg$ if there is no occurrence of $\sg$ in $\pi$.
One of the fundamental questions in the area of permutation patterns is to determine the number of permutations (words) of length $n$ containing $k$
occurrences of a given pattern $p$. That is, we want to find an
an explicit formula or a {\em generating function} for such permutations.
It is also of
great interest to find bijections between classes of
permutations and/or words that satisfy some sort of pattern condition
and other combinatorial
structures that preserve as many statistics as possible. For example,
we say say that two permutations $\sg, \tau \in S_k$ are {\em Wilf equivalent} if for each $n$,
the number of permutations in $S_n$ that avoid $\sg$ equals the number of permutations in $S_n$
that avoid $\tau$. If $\sg$ and $\tau$ are Wilf equivalent, then
it is natural to ask for a bijection between the set of
permutations of $S_n$ which avoid $\sg$ and the set of permutations
of $S_n$ which avoid $\tau$ which preserves as many classical
permutation statistics as possible.   Such
statistics-preserving bijections not only reveal structural
similarities between different combinatorial objects, but they often
also reveal previously unknown properties of the structures being
studied.

In \cite{BS} Babson and Steingr\'{\i}msson introduced {\em
generalized permutation patterns} ({\em GPs}) that allow for
the requirement that
two adjacent letters in a pattern must be adjacent in the
permutation. If we write, say $2\mn 31$, then we mean that if this
pattern occurs in a permutation $\pi$, then the letters in $\pi$
that correspond to $3$ and $1$ are adjacent. For example, the
permutation $\pi=516423$ has only one occurrence of the GP
$2\mn 31$, namely the subword 564, whereas the GP $2\mn 3\mn
1$ occurs, in addition, in the subwords 562 and 563. Note that a pattern containing
a dash between each pair of consecutive letters is a classical pattern.

The motivation for introducing these patterns in~\cite{BS} was
the study of {\em Mahonian statistics}. Many interesting results on GPs appear in the literature
(see~\cite{St} for a survey). In particular, \cite{C} provides
relations of generalized patterns to several well studied
combinatorial structures, such as set partitions, {\em Dyck paths},
{\em Motzkin paths} and involutions.
We refer to~\cite{KM2} for a survey over results on patterns
discussed so far.

Further generalizations and refinements of GPs have appeared in the
literature. For example, one can study occurrences of a pattern
$\sg$  in a permutation $\pi$ where one pays attention to the parity
of the elements in the subsequences of $\pi$ which are order
isomorphic to $\sg$. For instance, Kitaev and Remmel \cite{KR1}
studied descents (the GP 21) where one fixes the parity of exactly
one element of a descent pair. Explicit formulas for the
distribution of these (four) new patterns were provided. The new
patterns are shown in~\cite{KR1} to be connected to the {\em
Genocchi numbers}, the study of which goes back to Euler.
In~\cite{KR2}, Kitaev and Remmel generalized the results
of~\cite{KR1} to classify descents according to equivalence mod $k$
for $k \geq 3$ of one of the descent pairs. As a result of this
study, one obtains, in particular, remarkable binomial identities.
Liese~\cite{L,L1} studied enumerating descents where the difference
between descent pairs is a fixed value. More precisely, study of the
set $Des_k(\sigma)=\{i|\sigma_i-\sigma_{i+1}=k\}$ is done
in~\cite[Chpt 7]{L1}. Hall and Remmel~\cite{HR} further generalized
the studies in~\cite{KR1,KR2,L,L1}. The main focus of~\cite{HR} is
to study the distribution of descent pairs whose top $\sigma_i$ lies
in some fixed set $X$ and whose bottom $\sigma_{i+1}$ lies in some
fixed set $Y$ (such descents are called ``$(X,Y)$-descents'').
Explicit inclusion-exclusion type formulas are given for the number
of $n$-permutations with $k$ $(X,Y)$-descents.

A new type of permutation pattern condition was introduced by
Bousquet-M\'elou et al. in \cite{BCDK}.   In \cite{BCDK},
the authors considered restrictions on both the places and the values
where a pattern can occur. For example, they considered
pattern diagrams as pictured in Figure \ref{figure:little}.

\fig{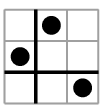}{A new $2~3~1$ pattern.}

\noindent The vertical line between the 2 and 3 of the pattern means
that the occurrence of the 2 and 3 in a subsequence must be
consecutive and the horizontal line between 2 and the 1 in the
pattern means that values corresponding to the 1 and 2 in an
occurrence must be consecutive. Thus an occurrence of the pattern
pictured above in a permutation $\pi = \pi_1 \ldots \pi_n \in S_n$,
is a sequence $1 \leq i_1 < i_2 < i_3 \leq n$ such that $\pi_{i_1}
\pi_{i_2} \pi_{i_3}$ is order isomorphic to $231$, $i_1+1 =i_2$ and
$\pi_{i_3} +1 = \pi_{i_1}$. For example, if $\pi = 31524$, then
there is an occurrence of $231$ in $\pi$, namely, $352$, but it is
not an occurrence of the pattern in Figure \ref{figure:little}
because 3 and 5 do not occur consecutively in $\pi$.  However if
$\tau = 32541$, then $251$ is an occurrence of the pattern as
pictured below.

\fig{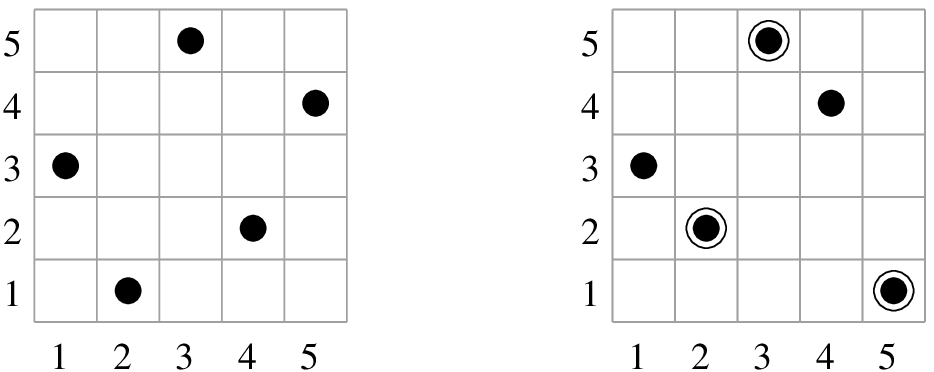}{Permutations which avoid and contain the new
$2~3~1$ pattern.}

\noindent
An attractive property of these
new patterns is that, like classical patterns (but not like GPs!), they
are closed under the action of $D_8$, the symmetry group of the
square. More precisely, the authors in~\cite{BCDK} studied
permutations that either avoid the pattern $2\mn 3\mn 1$ or in an
occurrence $\pi_i\pi_j\pi_k$ of $2\mn 3\mn 1$ in a permutation
$\pi_1\pi_2\ldots\pi_n$ where one either has $j\neq i+1$ or $\pi_i\neq\pi_k+1$.
It turns out that there is a bijection preserving several statistics
between {\em ({\bf 2}+{\bf 2})-free posets} and permutations
avoiding the pattern in previous sentence (see~\cite{BCDK}).

The outline of this paper is as follows. In Section~\ref{newconcept}
we define {\em place-difference-value patterns} ({\em PDVPs}) in
both permutations and words. These patterns cover under one roof
most of commonly used pattern restrictions that have occurred in the
literature on generalizations of GPs. In Sections~\ref{newres}
and~\ref{newresults} we consider several examples of PDVPs. Some of
these examples show connections to other combinatorial objects,
which cannot be obtained using the languages of most general notions
of GPs studied so far (we raise four bijective questions linking our
patterns to other combinatorial objects; see Problems 1--4).
Finally, in Section~\ref{final} we sketch some directions of further
research.

\section{Place-difference-value patterns}\label{newconcept}

In what follows, $\PP$ denotes the set of positive integers and
$k\PP$ denote the set of all positive multiples of $k$.

\begin{definition} A {\em place-difference-value pattern}, {\em PDVP}, is a
quadruple $P=(p,X,Y,Z)$ where $p$ is a permutation of length $m$,
$X$ is an $(m+1)$-tuple of non-empty, possibly infinite, sets of
positive integers, $Y$ is a set of triples $(s,t,Y_{s,t})$ where $0
\leq s < t \leq m+1$ and $Y_{s,t}$ is a non-empty, possibly
infinite, set of positive integers, and $Z$ is an $m$-tuple of
non-empty, possibly infinite, sets of positive integers. A PDVP
$P=(p_1p_2\ldots p_m,(X_0,X_1,\ldots,X_m),Y,(Z_1, \ldots, Z_m))$
occurs in a permutation $\pi=\pi_1\pi_2\ldots\pi_n$, if $\pi$ has a
subsequence $\pi_{i_1}\pi_{i_2}\ldots\pi_{i_m}$ with the following
properties:
\begin{enumerate}
\item $\pi_{i_k}<\pi_{i_{\ell}}$ if and only if $p_k < p_{\ell}$ for $1\leq k <\ell\leq m$;
\item $i_{k+1}-i_{k}\in X_k$ for $k=0,1,\ldots,m$, where we assume
$i_0=0$ and $i_{m+1}=n+1$;
\item for each $(s,t,Y_{s,t}) \in Y$, $|\pi_{i_s} - \pi_{i_t}| \in Y_{s,t}$ where we assume $\pi_{i_0} = \pi_0 =0$ and $\pi_{i_{m+1}} = \pi_{n+1} =n+1$; and
\item $\pi_{i_k} \in Z_k$ for $k =1,\ldots m$.
\end{enumerate}
\end{definition}

For example, let $\E$ and $\O$ denote the set of even and odd
numbers, respectively. Then the PDVP
$(12,(\{1\},\{3,4\},\{1,2,3\}),\{(1,2,\E)\},(\E,\PP))$ occurs in the
permutation $\pi=23154$ once as the subsequence 24. Indeed, each
such occurrence must start at position 1 as required by the set
$X_0$ and the second element of the sequence must occur either  at
position 4 or 5 as allowed by $X_1$. In this case $X_2$ does not
impose any additional restrictions. The condition that $Z_1 = \E$
says that the value in position $1$ must be even. Finally the
condition that $(1,2,\E) \in Y$ rules out the fact that the
subsequence $25$ is an occurrence of the pattern.

Classical
patterns are PDVPs of the form
$(p,(\PP,\PP,\ldots),\emptyset,(\PP,\PP,\ldots))$, whereas the GPs introduced
in~\cite{BS} have the property that $X_i$ is either $\PP$ or
$\{1\}$, $Y = \emptyset$, and $Z_i =\PP$ for all $i$. Also, the patterns introduced
in~\cite{BCDK} have the property that each of $X_i$'s
are either $\PP$ or $\{1\}$, $Z_i = \PP$ for all $i$, and
all the elements of $Y$ are of the form $(i,j,\{1\})$. Similarly,
the occurrences of the pattern $(21,(\PP,\{1\},\PP),\emptyset,(X,Y))$
in a permutation $\pi$ correspond to the $(X,Y)$-descents in
$\pi$ considered by
Hall and Remmel \cite{HR} and
the occurrences of the pattern $(21,(\PP,\{1\},\PP),\{(1,2,\{k\})\},(\PP,\PP))$
in $\pi$ correspond to elements of $Des_k(\pi)$ as studied by
Liese \cite{L,L1}.

We should note that there is often more than one way to specify the
same pattern. For instance, we can restrict ourselves to occurrences
of patterns that involve only even numbers by either setting $Z_i =
\E$ for all $i$ or by setting $Y = \{(i,i+1,\E):i=0, \ldots,m-1\}$.

In
Table~\ref{comparison}, we list how several pattern conditions
that have appeared in the literature can be expressed in terms
of PDVPs.

\begin{table}[htbp]
  \begin{center}
      \begin{tabular}{l|l}
        {\bf Object in the literature} & {\bf PDVP $P=(p,X,Y,Z)$}\\
        \hline
Classical patterns & $X_i=Z_j =\PP$ for all $i$ and $j$ and
$Y = \emptyset$.\\
\hline
GPs in~\cite{BS} & $X_i$ is either $\PP$ or $\{1\}$ for all $i$, $Z_j= \PP$ for all $j$, and $Y = \emptyset$.\\
\hline Descents conditioned on the
parity & $p =21$ and $X_0 = \PP$, $X_1 = \{1\}$, $X_2 =\PP$, \\
of the elements in the & $Y = \emptyset$, and $(Z_1,Z_2)$ equals \\
descent pairs as in \cite{KR1} & $(\E,\PP)$,  $(\O,\PP)$,
$(\PP,\E)$,or $(\PP,\O)$.\\
\hline Patterns in~\cite{KR2} & Similar to the last patterns,\\
& except we allow $Z_i$'s of the form $k \PP$ where $k \geq 3$.\\
\hline Patterns in~\cite{L,L1} & $(21,(\PP,\{1\},\PP),\{(1,2,\{k\})\},(\PP,\PP))$, where $k\geq 1$.\\
\hline Patterns in~\cite{HR} & $(21,(\PP,\{1\},\PP),\emptyset,(X,Y))$, $X$ and $Y$ are any fixed sets \\
\hline Patterns in~\cite{BCDK} & $X_i$ is either $\PP$ or $\{1\}$,
the elements of $Y$ are of the form \\
& $(i,j,\{1\})$ , and $Z_i = \PP$ for all $i$.\\
\hline
        \end{tabular}
  \end{center}
  \caption{Objects studied in the literature using the language of place-difference-value patterns.} \label{comparison}
\end{table}

The place-difference-value patterns in case of words can be defined in a
similar manner.

\begin{definition} A {\em place-difference-value (word) pattern}, {\em PDVP}, is a
quadruple
$P=(p,X,Y,Z)$ where $p$ is a word of length $m$ having an occurrence
of each of the letters $1,2,\ldots,k$ for some $k$, $X$ is an
$(m+1)$-tuple of non-empty, possibly infinite, sets of positive
integers, the elements of $Y$ are of the form
$(s,t,Y_{s,t})$ where $0 \leq i < j \leq m+1$ and
$Y_{i,j}$ is a non-empty, possibly infinite, set of non-negative integers, and  $Z$ is an $m$-tuple of non-empty, possibly infinite,
sets of positive
integers. A PDVP $P=(p_1p_2\ldots
p_m,(X_0,X_1,\ldots,X_m),Y,(Z_1, \ldots, Z_m))$ occurs in a word
$w=w_1w_2\ldots w_n$ over the alphabet $\{1,2,\ldots,t\}$, if $w$
has a subsequence $w_{i_1}w_{i_2}\ldots w_{i_m}$ with the following
properties:
\begin{enumerate}
\item $w_{i_k}<w_{i_{\ell}}$ if and only if $p_k<p_{\ell}$ for $1\leq k <\ell\leq
m$;
\item $w_{i_k}=w_{i_{\ell}}$ if and only if $p_k=p_{\ell}$ for $1\leq k <\ell\leq
m$;
\item $i_{k+1}-i_{k}\in X_k$ for $k=0,1,\ldots,m$, where we assume
$i_0=0$ and $i_{m+1}=n+1$;
\item for each $(s,t,Y_{s,t}) \in Y$, $|w_{i_s} - w_{i_t}| \in Y_{s,t}$ where we assume $w_{i_0} = w_0 =0$ and $w_{i_{m+1}} = w_{n+1} =t$; and
\item $w_{i_j} \in Z_j$ for $j =1, \ldots, m$.
\end{enumerate}
\end{definition}

We would like to point out that, of course, the notion  of PDVPs can
be generalized even further, e.g., increasing the number of dimensions
(as it is done in~\cite{KMV,KR}) or, for example,
 by having the differences or values
be dependent on the place. We will discuss some of these extensions
in Section~\ref{final}. In any case, the PDVPs are the closest
objects to those most popular pattern restrictions that have
appeared in the current literature on patterns.

Another thing to point out is that particular cases of patterns
introduced by us already appear in the literature, without any
general framework though. For example, in~\cite{T}, Tauraso found the
number of permutations of size $n$ avoiding simultaneously the
PDVPs
$$(12,(\PP,\{d\},\PP),\{(1,2,\{d\})\},(\PP,\PP))$$ and
$$(21,(\PP,\{d\},\PP),\{(1,2,\{d\})\},(\PP,\PP)),$$ where $2\leq d\leq n-1$.
Also, see~\cite[A110128]{oeis} for related objects.

\section{Some results on place-difference-value patterns in
permutations}\label{newres}

Recall that $\E = \{0, 2, 4, \ldots \}$ and $\O = \{1, 3, 5, \ldots
\}$ denote the set of even and odd numbers, respectively. Also, let
$\S_n$ denote the set of permutations of length $n$.

\subsection{Distribution of certain PDVPs on permutations}

Suppose $P=(p,(\O,\E, \ldots,\E), \emptyset, (\E, \ldots, \E))$, where $p$ is
any permutation on $t$ elements. We will show an easy connection
between distributions of $p$, viewed as a classical pattern, and the
pattern $P$.

The restriction on the $X$'s says that we want our pattern to occur
at the odd positions. The restriction on $Z$'s says that we are
worried about the even numbers that appear at the odd positions.

Let $A_{n,m}$ (resp. $B_{n,m}$) be the number of permutations in
$\S_n$ that contain $m$ occurrences of the pattern $p$ (resp. $P$).
There is an easy way to express $B_{n,m}$ in terms of $A_{n,m}$, as
it is shown below.

Consider $S_{2n}$ and suppose a permutation contains $m$ occurrences
of $P$  and $k$ even numbers in odd positions, where $0\leq k\leq
n$. We can choose the even numbers that appear in odd positions in
$\binom{n}{k}$ ways, and then we can choose the positions of the
even numbers in $\binom{n}{k}$ ways. Once we have chosen those
numbers and those positions, we have to arrange the even numbers so
that the permutation built by them would contain $m$ occurrences of
$p$. This can be done in $A_{k,m}$ ways. Then we have to choose the
odd numbers that appear in the odd positions in $\binom{n}{n-k}$
ways and those odd numbers can be arranged in $(n-k)!$ ways. Finally
the numbers which occupy the even positions can be arranged in $n!$
ways. Thus, \begin{equation}\label{result}B_{2n,m}=\sum_{k=0}^n
n!(n-k)!{n\choose k}^3A_{k,m}.\end{equation}
The case for $S_{2n+1}$ is similar.

Thus, whenever we know distribution of a classical pattern $p$ (it is
known only in a few cases), we can find distribution of $P$. On the
other hand, we can use the same formula for avoidance matters, in
which case we can get more applications of it. For example, for
$p=12$, $A_{k,0}$ becomes 1, as the only way to avoid $1\mn 2$ is to
arrange the corresponding even elements in decreasing order. Thus,
in this case, we have
$$B_{2n,0}=\sum_{k=0}^n n!(n-k)!{n\choose k}^3.$$
Another example is when $p=123$. It is well-known that the number of
$n$-permutations avoiding the pattern $1\mn2\mn 3$ is given by the
$n$-th Catalan number $C_n$, and thus, in this case, we have
$$B_{2n,0}=\sum_{k=0}^n n!(n-k)!{n\choose k}^3C_k.$$

For a slightly more complicated example, suppose that
$p =12$, $X = (\O,4\PP,\PP)$, $Y=\{(1,2,4\PP)\}$, and
$Z=(\O,\PP)$. Let $A = \{1,5,9,13, \ldots\}$ and
$B =\{3,7,11,15,\ldots\}$.  The restriction imposed by our choice of $X$ says
that we are only interested in subsequences  that occur at positions
in $A$  or subsequences that occur at positions in $B$.
 The restrictions imposed by our choice of
$Y$ and $Z$ says that we are only interested in subsequences that
involve values in $A$ or subsequences that involve values in $B$.
Suppose we want to find the number $K_{n}$ of permutations  in
$\S_{n}$ that avoid our pattern. Then choose $k_1$ to be the number
of elements of $A$ that occur in positions in $A$ and $k_2$ to be
the number of elements of $A$ that occur in positions in $B$.
Similarly choose $l_1$ to be the number of elements of $B$ that
occur in positions in $A$ and $l_2$ to be the number of elements of
$B$ that occur in positions in $B$.  To avoid our pattern, the $k_1$
elements of $A$ that occur in the positions of $A$ must be in
decreasing order and $k_2$ elements of $A$ that occur in the
positions of $B$ must occur in decreasing order. Next, the $l_1$
elements of $B$ that occur in the positions of $A$ must be in
decreasing order and $l_2$ elements of $B$ that occur in the
positions of $B$ must occur in decreasing order. Then we can arrange
the remaining elements that occur in the positions of $A$ in any
order we want and we can arrange in the remaining elements that
occur in the positions of $B$ any order we want.  Finally, we can
arrange the elements that lie in positions outside of $A$ and $B$ in
any order that we want. Thus, our final answer is determined by the
number of ways to choose the elements that correspond to $k_1$,
$k_2$, $l_1$ and $l_2$ and the number of ways to choose their
corresponding positions in $A$ and $B$. Thus, for example, one can
easily check that
\begin{eqnarray*}
&&K_{4n} = \\
&&(2n)!\sum_{\stackrel{0 \leq k_1,k_2,l_1,l_2}{k_1+k_2 \leq n, l_1+l_2 \leq n}} \binom{n}{k_1,k_2,n-k_1-k_2} \binom{n}{l_1,l_2,n-l_1-l_2} \times \\
&& \ \ \ \
\binom{n}{k_1,l_1,n-k_1-l_1}\binom{n}{k_2,l_2,n-k_2-l_2} (n-k_1-l_1)!
(n-k_2-l_2)! = \\
&& (2n)!\sum_{\stackrel{0 \leq k_1,k_2,l_1,l_2}{k_1+k_2 \leq n, l_1+l_2 \leq n}}\frac{(n!)^4}{(k_1!)^2(k_2!)^2(l_1!)^2(l_2!)^2(n-k_1-k_2)!(n-l_1-l_2)!}.
\end{eqnarray*}

\subsection{One more result on PDVPs on permutations}

In this subsection we consider the permutations which simultaneously avoid
 the GPs 231 and 132  and the PDVP
$P=(12,(\PP,\{k\},\PP),\{(1,2,\{1\})\},(\PP,\PP))$, where $k\geq 1$.
Let $a_{n,k}$ be the number of such permutations of length $n$. We
will show that
$$
a_{n,k} =
\begin{cases}
  F(n) & \text{~if~} k=1,\\
  2^{n-1} & \text{~if~} k\geq 2 \text{~and~} n\leq k,\\
  3\cdot 2^{n-3} & \text{~if~} k\geq 2 \text{~and~} n\geq k+1
\end{cases}
$$
where $F(n)$ is the $n$-th Fibonacci number.
The sequence of $a_{n,2}$'s --- 3, 6, 12, 24, 48, 96, 192,$\ldots$
--- appears in~\cite[A042950]{oeis}.

Notice, that avoiding just the GPs 231 and 132 gives $2^{n-1}$
permutations of length $n$ (\cite{K}), and the structure of such
permutations is a decreasing word followed by an increasing word (1
is staying between the words and it is assumed to belong to both of
them). Suppose first that $k\geq 2$. If $n\leq k$, then there is no
chance for $P$ to occur thus giving $2^{n-1}$ possibilities. On the
other hand, assuming $n=k+1$, the number of permutations avoiding
the three patterns is given by $2^{k-1}-2^{k-2}$ as whenever the
first letter is $n-1$ and the last letter is $n$, we get an
occurrence of $P$ (there are $2^{k-2}$ such cases). Finally, if we
increase the number of letters in a ``good'' permutations of length
$k+1$, one by one, we always have two places to insert a current
largest element: at the very beginning or at the very end, which
gives in total $(2^{k-1}-2^{k-2})2^{n-k-1}=3\cdot 2^{n-3}$
possibilities, as claimed.

In the case $k=1$, we think of counting good permutations by
starting with the letter 1, and inserting, one by one, the letters
2,3,$\ldots$. If $P$ would not be prohibited, we would always have
two choices to insert a current largest element. However, inserting
$n$, the configuration $(n-1)n$ is prohibited, which leads
immediately to a recursion for the Fibonacci numbers.

\begin{remark} Because of the structure of permutations avoiding GPs 231 and 132, one can see that the maximum number of occurrences of
$P$ in such permutations is 1. Thus, we actually found distribution
of $P$ on 231- and 132-avoiding permutations, as the number of such
permutations having exactly one occurrence of $P$ is
$2^{n-1}-a_{n,k}$.\end{remark}

An interpretation of the sequence \cite[A042950]{oeis}, based on a
result in~\cite{KMS}, suggests the following bijective question.

\begin{problem} For $k\geq 2$ and $n\geq k+1$, find a bijection between permutations of length $n$ which  simultaneously avoid
the GPs 231 and 132 and the PDVP
$P=(12,(\PP,\{k\},\PP),\{(1,2,\{1\})\},(\PP,\PP))$ and the set of
rises (occurrence of the GP 12) after $n$ iterations of the morphism
$1\rightarrow 123$, $2\rightarrow 13$, $3\rightarrow 2$, starting
with element 1. For example, for $k=2$ and $n=3$, there are 3
permutations avoiding the prohibitions, 123, 312, and 321, and there
are 3 rises in 1231323.
\end{problem}

\section{Some results on place-difference-value patterns on words}\label{newresults}

In this section, we consider examples of PDVPs on words  involving
both distance and value, that cannot be expressed in terms of
 pattern conditions that have appeared in the literature so far.

\subsection{The PDVP $(12,(\PP,\{2\},\PP),\{(1,2,\{2\})\},(\PP,\PP))$ on words.}
Consider words $w \in \{1, \ldots, k\}^*$. If $w=w_1 \ldots w_n$,
then let $S(w) =\{i:w_{i+2} -w_i =2\}$ and $s(w) = |S(w)|$.

Our goal is to compute the generating function
\begin{equation}\label{eq:kbasic}
A_k(q,z) = \sum_{w \in \{1, \ldots,k\}^*}q^{|w|}z^{s(w)}.
\end{equation}

Let
\begin{equation}\label{eq:2kbasic}
A_k(i_1 \ldots i_j;q,z) = \sum_{w \in \{1,
\ldots,k\}^*}q^{|i_1\ldots i_jw|}z^{s(i_1 \ldots i_jw)}.
\end{equation}
Then for example, when $k =3$, we easily obtain the following
recursions for $A_3(ij;q,z) = A(ij;q,z)$.
\begin{eqnarray*}
A(11;q,z) &=&  q^2 + qA(11;q,z) +qA(12;q,z) +qzA(13;q,z)\\
A(12;q,z) &=& q^2 + qA(21;q,z) +qA(22;q,z) +qzA(23;q,z)\\
A(13;q,z) &=& q^2 + qA(31;q,z) +qA(32;q,z) +qzA(33;q,z)\\
A(21;q,z) &=& q^2 + qA(11;q,z) +qA(12;q,z) +qA(13;q,z)\\
A(22;q,z) &=& q^2 + qA(21;q,z) +qA(22;q,z) +qA(23;q,z)\\
A(23;q,z) &=& q^2 + qA(31;q,z) +qA(32;q,z) +qA(33;q,z)\\
A(31;q,z) &=& q^2 + qA(11;q,z) +qA(12;q,z) +qA(13;q,z)\\
A(32;q,z) &=& q^2 + qA(21;q,z) +qA(22;q,z) +qA(23;q,z)\\
A(33;q,z) &=& q^2 + qA(31;q,z) +qA(32;q,z) +qA(33;q,z).
\end{eqnarray*}
Thus if we let
\begin{eqnarray*}
Q &=& (-q^2,-q^2,-q^2,-q^2,-q^2,-q^2,-q^2,-q^2,-q^2)\ \mbox{and}\\
A &=&
(A(11;q,z),A(12;q,z),A(13;q,z),A(21;q,z),A(22;q,z),A(23;q,z),\\
&& A(31;q,z), A(32;q,z),A(33;q,z)),
\end{eqnarray*}
then we see that
$$Q^T = MA^T$$
where
$$
M = \left( \begin{array}{ccccccccc}
q-1 & q  & zq  & 0  & 0  & 0  & 0 & 0 & 0 \\
0   & -1 & 0   & q  & q  & qz & 0 & 0 & 0 \\
0   & 0  & -1  & 0  & 0  & 0  & q & q & qz \\
q   & q  &  q  & -1 & 0  & 0  & 0 & 0 & 0  \\
0   & 0  & 0   & q  &q-1 & q  & 0 & 0 & 0 \\
0   & 0  & 0   & 0  & 0  & -1 & q & q & q\\
q   & q  & q   & 0  & 0  & 0  &-1 & 0 & 0 \\
0   & 0  & 0   & q  & q  & q  & 0 &-1 & 0\\
0   & 0  & 0   & 0  & 0  & 0  & q & q &q-1
\end{array} \right)
$$

Thus $A^T = M^{-1}Q^T$ and our desired generating function is given
by
\begin{equation}
A_3(q,z) = 1 + 3q + (1,1,1,1,1,1,1,1,1)A^T = \frac{1}{(1-q^2(1-z))(1
-(3q + q^2(z-1)))}.
\end{equation}
where we used Mathematica for the last equation.

Note that
\begin{eqnarray*}
\frac{1}{1 -(3q + q^2(z-1))} &=&  \sum_{m \geq 0}  \sum_{k=0}^m
\binom{m}{k}3^kq^kq^{2m-2k}(z-1)^{m-k} \\
&=&  \sum_{m \geq 0}  \sum_{k=0}^m
\binom{m}{k}3^k(z-1)^{m-k}q^{2m-k}.
\end{eqnarray*}
Now it is easy to see that $\sum_{k=0}^m
\binom{m}{k}3^k(z-1)^{m-k}q^{2m-k}$ involves powers of $q$ that
range from $m$ to $2m$. Thus, since $2m-k =2n$ if and only if $k =
2m -2n$, we have
\begin{eqnarray*}
\frac{1}{1 -(3q + q^2(z-1))}|_{q^{2n}} &=&
\sum_{m=n}^{2n} \binom{m}{2m-2n} 3^{2m-2n}(z-1)^{m-(2m-2n)} \\
&=& \sum_{m=0}^n \binom{m+n}{2m} 3^{2m} (z-1)^{n-m}.
\end{eqnarray*}
It follows that
\begin{eqnarray*}
\frac{1}{(1-q^2(1-z))(1 -(3q + q^2(z-1))}|_{q^{2n}} &=&
\sum_{r=0}^n (1-z)^{n-r} \sum_{m=0}^r \binom{m+r}{2m} 3^{2r} (z-1)^{r-m} \\
&=& \sum_{r=0}^n  \sum_{m=0}^r (-1)^{n-r} \binom{m+r}{2m} 3^{2m}
(z-1)^{n-m}
\end{eqnarray*}
and
\begin{equation}
\frac{1}{(1-q^2(1-z))(1 -(3q + q^2(z-1))}|_{q^{2n}z^s} =
\sum_{r=0}^n  \sum_{m=0}^r (-1)^{m+r+s} \binom{m+r}{2m}
\binom{n-m}{s}9^m.
\end{equation}
Similarly, one can show that
\begin{eqnarray*}
\frac{1}{1 -(3q + q^2(z-1))}|_{q^{2n+1}} &=&
\sum_{m=n+1}^{2n+1} \binom{m}{2m-2n-1} 3^{2m-2n-1}(z-1)^{m-(2m-2n-1)} \\
&=& \sum_{m=0}^n \binom{m+n+1}{2m+1} 3^{2m+1} (z-1)^{n-m}.
\end{eqnarray*}
It follows that
\begin{eqnarray*}
\frac{1}{(1-q^2(1-z))(1 -(3q + q^2(z-1))}|_{q^{2n+1}} &=&
\sum_{r=0}^n (1-z)^{n-r} \sum_{m=0}^r \binom{m+r+1}{2m+1} 3^{2r+1}(z-1)^{r-m} \\
&=& \sum_{r=0}^n  \sum_{m=0}^r (-1)^{n-r} \binom{m+r+1}{2m+1}
3^{2m+1} (z-1)^{n-m}
\end{eqnarray*}
and
\begin{equation}
\frac{1}{(1-q^2(1-z))(1 -(3q + q^2(z-1))}|_{q^{2n+1}z^s} =
\sum_{r=0}^n  \sum_{m=0}^r (-1)^{m+r+s} \binom{m+r+1}{2m+1}
\binom{n-m}{s}3^{2m+1}.
\end{equation}

Thus we have shown that, in particular, in case of avoidance,
\begin{eqnarray}
A_3(q,0)|_{q^{2n}} &=&  \sum_{r=0}^n  \sum_{m=0}^r (-1)^{m+r}
\binom{m+r}{2m}
3^{2m} \ \mbox{and}\\
A_3(q,0)|_{q^{2n+1}} &=&  \sum_{r=0}^n  \sum_{m=0}^r (-1)^{m+r}
\binom{m+r+1}{2m+1} 3^{2m+1}.
\end{eqnarray}

A check in~\cite{oeis} shows that $A_3(q,0)|_{q^{2n}} =(F(2n))^2$
and
 $A_3(q,0)|_{q^{2n+1}} =F(2n)F(2n+2)$ where $F(n)$ is the {\em $n$-th Fibonacci number} (see~\cite[A006190]{oeis}). Here is a proof of that fact.
 Clearly if we have a word
$w=w_1 \ldots w_{2n}$ such that $s(w) =0$, then $u = w_1 w_3 \ldots
w_{2n-1}$ and $v = w_2 w_4 \ldots w_{2n}$ must be words in
$\{1,2,3\}^*$ that never have a 3 following a 1. The map
$\{1\rightarrow 01, 2\rightarrow 00, 3\rightarrow 10\}$ gives a
bijection from the set of words of length $n$ avoiding 13 and the
set of binary words of length $2n$ avoiding 11 and known to be
counted by $F(2n)$.

For another way to understand the same result,
we first find the distribution of consecutive occurrences of
13 over words in $\{1,2,3\}^*$. For any word $u
= u_1 \ldots u_n \in \{1,2,3\}^*$, let
$T(w) = \{i:w_{i+1}=2+w_i\}$ and $t(w) = |T(w)|$. Then
we wish to compute
\begin{equation}\label{eq:3basic}
B_3(q,z) = \sum_{w \in \{1,2,3\}^*}q^{|w|}z^{t(w)}.
\end{equation}
Let
\begin{equation}\label{eq:23basic}
B_3(i_1 \ldots i_j;q,z) = \sum_{w \in \{1,2,3\}^*}q^{|i_1\ldots
i_jw|}z^{t(i_1 \ldots i_jw)}.
\end{equation}
Then it is easy to see that
\begin{eqnarray*}
B_3(1;q,z) &=& q + qB_3(1;q,z) + qB_3(2;q,z) +qzB_3(3;q,z)\\
B_3(2;q,z) &=& q + qB_3(1;q,z) + qB_3(2;q,z) +qB_3(3;q,z)\\
B_3(3;q,z) &=& q + qB_3(1;q,z) + qB_3(2;q,z) +qB_3(1;q,z).\\
\end{eqnarray*}
Thus if $\bar{Q} = (-q,-q,-q)$ and $B =
(B_3(1;q,z),B_3(2;q,z),B_3(3;q,z))$, then $\bar{Q}^T = R B^T$ where
$$
R = \left( \begin{array}{ccc} q-1 & q & qz \\ q & q-1 & q \\ q & q &
q-1
\end{array} \right).
$$
Thus $B^T =R^{-1}\bar{Q}^T$ and hence
\begin{eqnarray}\label{B3}
B_3(q,z) &=&  1 + B_3(1;q,z) + B_3(2;q,z) +  B_3(3;q,z) \\
&=&  \frac{1}{1-3q -q^2(z-1)}.
\end{eqnarray}

To derive the avoidance case algebraically, notice that the
generating function for the Fibonacci numbers (with a proper shift
of indices) is
\begin{equation}\label{Fib}
F(q) = \sum_{n \geq 0} F(n)q^n = \frac{1+q}{1-q -q^2}.
\end{equation}
Thus
\begin{eqnarray*}
f(q) &=& \sum_{n \geq 0} F(2n) q^n \\
&=& \frac{F(q^{1/2}) + F(-q^{1/2})}{2} \\
&=& \frac{1}{1-3q +q^2} \\
&=& B_3(q,0).
\end{eqnarray*}
This shows again that the number of words $w \in \{1,2,3\}^*$ of
length $n$ such that  $t(w) =0$ is equal to $F(2n)$. It easily
follows that $w \in \{1,2,3\}^*$ of length $2n$ such that  $s(w) =0$
is equal to $(F(2n))^2$ and the number of word $w \in \{1,2,3\}^*$
of length $2n+1$ such that  $s(w) =0$ is equal to $F(2n)F(2n+2)$.

%\begin{equation}\label{Fib}
%F(q) = \sum_{n \geq 0} F(n)q^n = \frac{q}{1-q -q^2}.
%\end{equation}
%Thus
%\begin{eqnarray*}
%f(q) &=& \sum_{n \geq 0} F(2n) q^n \\
%&=& \frac{F(q^{1/2}) + F(-q^{1/2})}{2} \\
%&=& \frac{1}{1-3q +q^2} \\
%&=& B_3(q,0).
%\end{eqnarray*}
%This shows that the number of words $w \in \{1,2,3\}^*$ of length
%$n$ such that  $t(w) =0$ is equal to $F(2n)$. It easily follows that
%$w \in \{1,2,3\}^*$ of length $2n$ such that  $s(w) =0$ is equal to
%$(F(2n))^2$ and the number of word $w \in \{1,2,3\}^*$ of length
%$2n+1$ such that  $s(w) =0$ is equal to $F(2n)F(2n+2)$.

One can do a similar calculations when $k=4$. In that case,
Mathematica shows that
\begin{equation}
A_4(q,z) = \frac{1}{1-4q-8q^3(z-1)-4q^4(z-1)^2}
\end{equation}
and thus,
$$A_4(q,0) = 1+4q+16q^2+56q^3+196q^4+672q^5+2304q^6+ \cdots .$$
As before, if we have a word of $w=w_1 \ldots w_{2n}$ such that
$s(w) =0$, then $u = w_1 w_3 \ldots w_{2n-1}$ and $v = w_2 w_4
\ldots w_{2n}$ must be words in $\{1,2,3,4\}^*$ that never have a 3
following a 1 or a 4 following a 2. For any word $u = u_1 \ldots u_n
\in \{1,2,3,4\}^*$, let, as before, $T(w) = \{i:w_{i+1}=2+w_i\}$ and
$t(w) = |T(w)|$.  Then we can compute
\begin{equation}\label{eq:4basic}
B_4(q,z) = \sum_{w \in \{1,2,3,4\}^*}q^{|w|}z^{t(w)}
\end{equation}
in the same way that we computed $B_3(q,z)$. In this case,
\begin{equation}
B_4(q,z) = \frac{1}{1-4q -2q^2(z-1)}.
\end{equation}
Then
$$B_4(q,0)= 1+4q+14q^2+48q^3+164q^4+560q^5+1912q^6 +  \cdots.$$

It should be noted that $B_4(q,0) = \frac{1}{1-4q+2q^2}$ is a
generating function that beyond the objects listed
in~\cite[A007070]{oeis} (e.g., the number of words $w_0w_1\ldots
w_{2n+3}$ over $\{1,2,\ldots,7\}^*$ with $w_0=1$ and $w_{2n+3}=4$
and $|w_i-w_{i-1}|=1$) counts the number of independent sets in
certain ``almost regular'' graphs $G_3^n$ (see~\cite{BKM}). We leave
establishing a bijection between the objects in question as an open
problem, stating explicitly the one related to another open question
below.

\begin{problem} Find a bijection between the set $A_{n}$ of words $w=w_1w_2\ldots w_n\in \{1,2,3,4\}^*$ that avoid the pattern
$(12,(\PP,\{1\},\PP),\{(1,2,\{2\})\},(\PP,\PP))$ and the set $B_{n}$
of words $w_0w_1\ldots w_{2n+3}$ over $\{1,2,\ldots,7\}^*$ with
$w_0=1$ and $w_{2n+3}=4$ and $|w_i-w_{i-1}|=1$.
\end{problem}

\subsection{The PDVP $(12,(\PP,\{1,2\},\PP),\{(1,2,\{2\})\},(\PP,\PP))$ on words.}

Let $U(w) =\{i:w_{i+1} -w_i =2\}$ and $V(w) =\{i:w_{i+2} -w_i =2\}$
and let $p(w) = |U(w)|+|(V(w)|$. In that case, we can use
essentially the same methods to calculate $D_k(q,z)= \sum_{w \in
\{1,\ldots,k\}^*}q^{|w|}z^{p(w)}$.

Let
\begin{equation}\label{eq:2kmodified}
D_k(i_1 \ldots i_j;q,z) = \sum_{w \in \{1,
\ldots,k\}^*}q^{|i_1\ldots i_jw|} z^{p(i_1 \ldots i_jw)}.
\end{equation}
Then for example, when $k =3$, we easily obtain the following
recursions for $D_3(ij;q,z) =D(ij;q,z)$.
\begin{eqnarray*}
D(11;q,z) &=&  q^2 + qD(11;q,z) +qD(12;q,z) +qzD(13;q,z)\\
D(12;q,z) &=& q^2 + qD(21;q,z) +qD(22;q,z) +qzD(23;q,z)\\
D(13;q,z) &=& zq^2 + zqD(31;q,z) +zqD(32;q,z) +qz^2D(33;q,z)\\
D(21;q,z) &=& q^2 + qD(11;q,z) +qD(12;q,z) +qD(13;q,z)\\
D(22;q,z) &=& q^2 + qD(21;q,z) +qD(22;q,z) +qD(23;q,z)\\
D(23;q,z) &=& q^2 + qD(31;q,z) +qD(32;q,z) +qD(33;q,z)\\
D(31;q,z) &=& q^2 + qD(11;q,z) +qD(12;q,z) +qD(13;q,z)\\
D(32;q,z) &=& q^2 + qD(21;q,z) +qD(22;q,z) +qD(23;q,z)\\
D(33;q,z) &=& q^2 + qD(31;q,z) +qD(32;q,z) +qD(33;q,z).
\end{eqnarray*}
Thus if we let
\begin{eqnarray*}
Q &=& (-q^2,-q^2,-zq^2,-q^2,-q^2,-q^2,-q^2,-q^2,-q^2)\ \mbox{and}\\
D &=&
(D(11;q,z),D(12;q,z),D(13;q,z),D(21;q,z),D(22;q,z),D(23;q,z),\\
&& D(31;q,z),D(32;q,z),D(33;q,z)),
\end{eqnarray*}
then we see that
$$Q^T = MD^T$$
where
$$
M = \left( \begin{array}{ccccccccc}
q-1 & q  & zq  & 0  & 0  & 0  & 0 & 0 & 0 \\
0   & -1 & 0   & q  & q  & qz & 0 & 0 & 0 \\
0   & 0  & -1  & 0  & 0  & 0  & zq & zq & qz^2 \\
q   & q  &  q  & -1 & 0  & 0  & 0 & 0 & 0  \\
0   & 0  & 0   & q  &q-1 & q  & 0 & 0 & 0 \\
0   & 0  & 0   & 0  & 0  & -1 & q & q & q\\
q   & q  & q   & 0  & 0  & 0  &-1 & 0 & 0 \\
0   & 0  & 0   & q  & q  & q  & 0 &-1 & 0\\
0   & 0  & 0   & 0  & 0  & 0  & q & q &q-1
\end{array} \right)
$$

Thus $D^T = M^{-1}Q^T$ and our desired generating function is given
by
\begin{eqnarray}
D_3(q,z) &=& 1 + 3q + (1,1,1,1,1,1,1,1,1)D^T \nonumber \\
&=&
\frac{1}{1-3q-q^2(z-1)-q^3(2z+1)(z-1) -q^4(z-1)^2}
\end{eqnarray}
where we again used Mathematica for the last equation. In this case,
$$D_3(q,0)= 1+3q+8q^2+20q^3+49q^4+119q^5+288q^6+ \cdots.$$
The sequence $1,3,8,20,49,119,288, \ldots$ appears
in~\cite[A048739]{oeis} where it is given an interpretation as the
number of words $w=w_0 \ldots w_{n+2} \in \{1,2,3\}^*$ such that
$|w_{i+1}-w_i| \leq 1$, and $w_0 =1$ and $w_{n+2} =3$. We leave the
following as an open question:

\begin{problem} Find a bijection between the set $A_{n}$ of words $w=w_1w_2\ldots w_n\in \{1,2,3\}^*$ that avoid the pattern
$(12,(\PP,\{1,2\},\PP),\{(1,2,\{2\})\},(\PP,\PP))$ and the set
$B_{n}$ of words $w=w_0 \ldots w_{n+1} \in \{1,2,3\}^*$ such that
$|w_{i+1}-w_i| \leq 1$, and $w_0 =1$ and $w_{n+1} =3$. \end{problem}

%\begin{proof}
%%We first show that both of the sequences $|A_n|$ and $|B_n|$ satisfy the same
%%recursion
%%$$a_n=2a_{n-1}+a_{n-2}+1$$ where $a_n=|A_n|=|B_{n+1}|$. Also,
%%$a_1=3$. Indeed, for $A_n$, either a word begins with 2 or 3 which
%%does not bring any restriction on the remaining $n-1$ letters (and
%%thus gives the term $2a_{n-1}$, or it begins with 1.
%
%We proceed with mapping $f:A_n\rightarrow B_{n+1}$ by induction on
%$n$ with the base case $1 \rightarrow 112$, $2 \rightarrow 121$, and
%$3 \rightarrow 122$. Let $w=w_1w_2\ldots w_n\in A_n$ and $n\geq 3$.
%We distinguish four non-overlapping cases that cover all the
%possibilities:
%\begin{itemize}
%\item $w_1=1$, and thus $w_2\neq 3$ and $w_3\neq 3$: $f(w)=11w_2w_3
%f'(w_4w_5\ldots w_n)$ where $f'(X)$ removes the leftmost letter in
%$f(X)$, and by definition, $f'(\epsilon)=3$ ($\epsilon$ is the empty
%word). Note that $f(w)\in B_n$ as $11w_2w_3$ is a proper factor, as
%well as $f'(w_4w_5\ldots w_n)$ (by induction), and the fact that
%$w_3\neq 3$ and the leftmost letter of $f'(w_4w_5\ldots w_n)$ cannot
%be 3, we have no problems with adjoining the factors to each other.
%\item $w_1=3$: $f(w)=122w_2f'(w_4w_5\ldots w_n)$.
%\item $f(2311\ldots 1)=1211\ldots 1$.
%\item $w_1=2$, $w_2=3$: $f(w)=12f(w_3w_4\ldots w_n)$.
%\end{itemize}
%\end{proof}

A similar computation as that above will show that
\begin{equation}
D_4(q,z) =
\frac{1+2q^2(1-z)-2q^3(z-1)^2}{1-4q-8q^2(z-1)-4q^4(z-1)^2}.
\end{equation}
It follows that
\begin{eqnarray*}
D_4(q,0)&=& \frac{1+2q^2-2q^3}{1-4q+8q^3-4q^4} \\
&=&  1+4q+14q^2+46q^3+156x^4+528x^5+1800x^6+ \cdots.
\end{eqnarray*}

\subsection{The PDVP $(12,(\PP,\{1,2\},\PP),\{(1,2,\{2\})\},(\O,\PP))$ on words.}

In this case, an occurrence of our PDVP is either 2 consecutive
odd numbers that differ by 2 or two odd numbers at distance 2 from
each other that differ by 2.
Let $P(w) =\{i:w_{i+1} -w_i =2 \ \& \ w_{i+1} \in \O\}$ and
$Q(w) =\{i:w_{i+2} -w_i =2 \ \& \ w_{i+2} \in \O\}$
and let $r(w) = |P(w)|+|(Q(w)|$. In that case, we can use
essentially the same methods to calculate $E_k(q,z)= \sum_{w \in
\{1,\ldots,k\}^*}q^{|w|}z^{r(w)}$.

Let
\begin{equation}\label{eq:2Ekmodified}
E_k(i_1 \ldots i_j;q,z) = \sum_{w \in \{1,
\ldots,k\}^*}q^{|i_1\ldots i_jw|} z^{r(i_1 \ldots i_jw)}.
\end{equation}

For example, in the case that $n=4$, it is easy to see that
\begin{eqnarray*}
E(11;q,z) &=&  q^2 + qE(11;q,z) +qE(12;q,z) +qzE(13;q,z)+qE(14:q,z)\\
E(12;q,z) &=&  q^2 + qE(21;q,z) +qE(22;q,z) +qzE(23;q,z)+qE(24:q,z)\\
E(13;q,z) &=&  q^2z + qzE(31;q,z) +qzE(32;q,z) +qz^2E(33;q,z)+qzE(34:q,z)\\
E(14;q,z) &=&  q^2 + qE(41;q,z) +qE(42;q,z) +qzE(13;q,z)+qE(14:q,z)\\
\end{eqnarray*}
and for any $i \in {2,3,4}$,
\begin{eqnarray*}
E(i1;q,z) &=&  q^2 + qE(11;q,z) +qE(12;q,z) +qE(13;q,z)+qE(14:q,z)\\
E(i2;q,z) &=&  q^2 + qE(21;q,z) +qE(22;q,z) +qE(23;q,z)+qE(24:q,z)\\
E(i3;q,z) &=&  q^2 + qE(31;q,z) +qE(32;q,z) +qE(33;q,z)+qE(34:q,z)\\
E(i4;q,z) &=&  q^2 + qE(41;q,z) +qE(42;q,z) +qE(13;q,z)+qE(14:q,z).\\
\end{eqnarray*}

Thus if we let
\begin{eqnarray*}
Q &=& (-q^2,-q^2,-zq^2,-q^2,-q^2,-q^2,-q^2,-q^2,-q^2,-q^2,-q^2,-q^2,-q^2,-q^2,-q^2,-q^2)\ \mbox{and}\\
E &=& (E(ij;q,z))_{1 \leq i,j \leq 4}
\end{eqnarray*}
where we order the elements of $E$ according the lexicographic order
on the pairs $(i,j)$, then we see that
$$Q^T = ME^T$$
where

{\tiny
$$
M = \left( \begin{array}{ccccccccccccccccc}
q-1 & q  & zq  & q  & 0 & 0  & 0  & 0  & 0  & 0  & 0  & 0 & 0  & 0  & 0  & 0 \\
0   & -1 & 0   & 0  & q & q  & qz & q  & 0  & 0  & 0  & 0 & 0  & 0  & 0  & 0 \\
0   & 0  & -1  & 0  & 0 & 0  & 0  & 0  & qz  & qz  & qz^2 & qz & 0  & 0  & 0  & 0 \\
0   & 0  & 0   & -1 & 0 & 0  & 0  & 0  & 0  & 0  & 0  & 0 & q  & q  &qz  & q \\
q   & q  & q   & q  &-1 & 0  & 0  & 0  & 0  & 0  & 0  & 0 & 0  & 0  & 0  & 0 \\
0   & 0  & 0   & 0  & q & q-1& q  & q  & 0  & 0  & 0  & 0 & 0  & 0  & 0  & 0 \\
0   & 0  & 0   & 0  & 0 & 0  & -1 & 0  & q  & q  & q  & q & 0  & 0  & 0  & 0 \\
0   & 0  & 0   & 0  & 0 & 0  & 0  & -1& 0  & 0  & 0  & 0 & q  & q  & q  & q \\
q   & q  & q   & q  & 0 & 0  & 0  & 0 & -1 & 0  & 0  & 0  & 0  & 0  & 0  & 0 \\
0   & 0  & 0   & 0  & q & q  & q  & q & 0  & -1 & 0  & 0  & 0  & 0  & 0  & 0 \\
0   & 0  & 0   & 0  & 0 & 0  & 0  & 0 & q  & q  &q-1 & q  & 0  & 0  & 0  & 0 \\
0   & 0  & 0   & 0  & 0 & 0  & 0  & 0 & 0  & 0  & 0  &-1  & q  & q  & q  & q \\
q   & q  & q   & q  & 0 & 0  & 0  & 0 & 0  & 0  & 0  & 0 & -1 & 0  & 0  & 0 \\
0   & 0  & 0   & 0  & q & q  & q  & q & 0  & 0  & 0  & 0 & 0  & -1 & 0  & 0 \\
0   & 0  & 0   & 0  & 0 & 0  & 0  & 0 & q  & q  & q  & q & 0  & 0  & -1 & 0 \\
0   & 0  & 0   & 0 & 0  & 0  & 0  & 0 & 0  & 0  & 0  & 0 & q  & q  & q  & q-1 \\\end{array} \right)
$$}
Thus $E^T = M^{-1}Q^T$ and our desired generating function is given
by
\begin{eqnarray}
E_4(q,z) &=& 1 + 4q + (1,1,1,1,1,1,1,1,1,1,1,1,1,1,1,1)E^T \nonumber \\
&=&
\frac{1}{1-4q -(z-1)q^2 - 2(z^2-1)q^3-z(z-1)^2 q^4 }
\end{eqnarray}
where we again used Mathematica for the last equation.
Then
$$E_4(q,0) = 1 +4q +15q^2 + 54q^3 +193q^4+ 688q^5 + \cdots .$$

\subsection{One more result on PDVPs on words}

In this subsection, we consider the following two PDVP's,
\begin{eqnarray*}
P_1 &=& (12,(\PP,\{1\},\PP),\{(1,2,\{1\})\},(\PP,\PP)) \
\mbox{and}\\
P_2 &=& (12,(\PP,\{2\},\PP),\{(1,2,\{2\})\},(\PP,\PP)).
\end{eqnarray*}
Our goal is to  show that the number $a_{n}$ of words of
length $n$ over $\{1,2,3\}$ avoiding simultaneously the PDVP
$P_1$ and
$P_2$ is given by
$F(n+4)-n-3$, where $F(n)$, as above, is the $n$-th Fibonacci
number. The corresponding sequence --- 3, 7, 14, 26, 46, 79, 133,
221, $\ldots$
--- appears in~\cite[A079921]{oeis}.

It is easy to check that $a_1=3$ and $a_2=7$. We will show that for
$n\geq 3$, $a_n = a_{n-1}+a_{n-2}+n+1$, thus proving the claim
by~\cite[A079921]{oeis}. For a given word $w=w_1w_2\ldots w_n$ over
$\{1,2,3\}$ avoiding the prohibited patterns, we distinguish 5
non-overlapping cases covering all possibilities:
\begin{enumerate}
\item\label{case1} $w_1=3$. This 3 has no effect on the rest of $w$ (it cannot be involved in an occurrence of a prohibited pattern) thus giving
$a_{n-1}$ possibilities.

\item $w_1w_2=11$. There is only one valid extension of 11
to the right, namely, $1\ldots 1$, as placing 2 (resp. 3) will
introduce an occurrence of $P_1$ (resp. $P_2$) in $w$. Thus, the
number of possibilities in this case is 1.

\item $w_1w_2=13$. Extending 13 to the right by any legal word $w_3w_4\ldots w_n$ of
length $n-2$, we will be getting valid words of length $n$, except
for the case when $w_3=3$ ($w_1w_3$ is an occurrence of $P_2$). The
number of ``bad'' words, according to case~(\ref{case1}) above is
$a_{n-3}$. Thus, the number of possibilities in this case is
$a_{n-2}-a_{n-3}$.

\item $w=\underbrace{2\ldots 2}_{>0}\underbrace{1\ldots 1}_{\geq
0}$. The number of possibilities in this case is, clearly, $n$.

\item\label{case5} $w_1=2$ and $w$ contains at least one 3. Notice, that the leftmost 3 in $w$ must be preceded by 1 (to avoid $P_1$),
which, in turn, must be preceded by 2 (using the fact that $w$
avoids $P_2$ and $w_1=2$). Thus, in this case, $w$ begins with a
word of the form $\underbrace{2\ldots 2}_{>0}13x$, where $x$, if it
exists, is not equal to 3. To count all such words, we proceed according
to the
following, obviously reversible, procedure. Consider a good word,
say $v$, of length $n-3$. If $v$ does not begin with 3, map it to
$213v$ to get a proper word of length $n$ in the class in question.
On the other hand, if $v=\underbrace{3\ldots 3}_{i>0}xV$, were
$x\neq 3$ (assuming such $x$ exists), map $v$ to
$\underbrace{2\ldots 2}_{i>1}13xV$ getting a proper word of length
$n$ in the class. Clearly, we get all words in the class. Thus, the
number of possibilities in this case is $a_{n-3}$.
\end{enumerate}

Summarizing cases (\ref{case1})--(\ref{case5}) above we get the
desired.

The problem below involves so called {\em 2-stack sortable
permutations}, that is, permutations that can be sorted by passing
them twice through a {\em stack} (where the letters on the stack
must be in increasing order). Such permutations were first
considered in~\cite{W}, but attracted much of attention in the
literature since then.

\begin{problem} Find a bijection between the set $A_{n}$ of words $w=w_1w_2\ldots w_n\in \{1,2,3\}^*$ that avoid simultaneously the
PDVPs $P_1$ and
$P_2$ and the set $B_{n}$
of 2-stack sortable permutations which avoid the pattern $1\mn3\mn
2$ and contain exactly one occurrence of the pattern $1\mn 2\mn 3$.
The last object is studied in~\cite{EM}.
\end{problem}

\section{Beyond PDVPs: directions of further research}\label{final}

    Another generalization of the GPs is {\em partially ordered patterns} ({\em POPs}) when the letters of
    a pattern form a partially ordered set (poset), and an occurrence of such a pattern in a permutation is
    a linear extension of the corresponding poset in the order suggested by the pattern (we also pay attention
    to eventual dashes and brackets). For instance, if we have a poset on three elements labeled by $1^{\prime}$, $1$, and $2$,
    in which the only relation is $1<2$ (see figure~\ref{poset01}), then in an occurrence of $p=1^{\prime}\mn 12$ in
    a permutation $\pi$ the letter corresponding to the $1^{\prime}$ in $p$ can be either larger or smaller than
    the letters corresponding to $12$. Thus, the permutation 31254 has three occurrences of $p$, namely $3\mn12$,
    $3\mn25$, and $1\mn25$.

\setlength{\unitlength}{4mm}
    \begin{figure}[h]
\begin{center}
\begin{picture}(6,2)
\put(0,0){\put(2,0){\p} \put(4,0){\p} \put(4,2){\p} \path(4,0)(4,2)

\put(1.2,0){$1'$} \put(4.5,0){1} \put(4.5,2){2}
 }
\end{picture}
\caption{A poset on three elements with the only relation $1<2$.}
\label{poset01}
\end{center}
\end{figure}
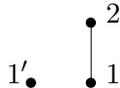

The notion of a POP allows to collect under one roof (to provide a
uniform notation for) several combinatorial structures such as {\em
peaks}, {\em valleys}, {\em modified maxima} and {\em modified
minima} in permutations, {\em Horse permutations} and $p$-{\em
descents} in permutations. See \cite{K1,K2,K3} for results,
including a survey paper, on POPs in permutations and \cite{KM1} on
POPs on words.

In the literature on permutation patterns, there are several
publications involving so called {\em barred patterns}. For example,
in~\cite{BCDK} a conjecture is settled on the number of permutations
avoiding the barred pattern $3\mn\bar{1}\mn5\mn2\mn\bar{4}$. A
permutation $\pi$ avoids $3\mn\bar{1}\mn5\mn2\mn\bar{4}$ if every
occurrence of the pattern $2\mn 3\mn1$ plays the role of 352 in an
occurrence of the pattern $3\mn1\mn5\mn2\mn4$. In some cases, barred
patterns can be expressed in terms of generalized patterns. E.g., to
avoid $4\mn1\mn\bar{3}\mn5\mn2$ is the same as to avoid $3\mn14\mn
2$. However, in many cases, one cannot express the barred patterns
in terms of other patterns. The pattern
$3\mn\bar{1}\mn5\mn2\mn\bar{4}$ is an example of such pattern.
Another example is the barred pattern $3\mn\bar{5}\mn2\mn4\mn1$ (it
is shown in~\cite{W} that the set of 2-stack sortable permutations
mentioned above is described by avoidance of
$3\mn\bar{5}\mn2\mn4\mn1$ and $2\mn 3\mn4\mn1$). In general, one can
consider distributions, rather than just avoidance, of barred
patterns. For example, the pattern $2\mn\bar{3}\mn1$ occurs in a
permutation $\pi$ $k$ times, if there are exactly $k$ occurrences
$ba$ in $\pi$ of the pattern $2\mn 1$ such that there is no element
$c>b$ in $\pi$ between $b$ and $a$.

It is straightforward to define {\em place-difference-value
partially order patterns}, {\em PDVPOPs}, or {\em
place-difference-value barred patterns}, {\em PDVBP}, since our
place, difference, and value restrictions just limit where and what
values are required for a pattern match. In particular,
formula~(\ref{result}) holds for PDVPOPs. We shall not pursue the
study of PDVPOPs or PDVBPs in this paper. Instead, we shall leave it
as a topic for further research.

Finally, we should observe that our definition of PDVP's does not
cover  the most general types of restrictions on patterns that one
can consider. For example, one can easily imagine cases where there
are restrictions on the values in occurrences of patterns that are a
function of the places occupied by the  occurrence or there are
restrictions on the places which an occurrence occupied that are
function  of the values in the occurrence. Thus the most general
type of restriction for a pattern $p \in S_m$ would be to just give
a set $\mathcal{S}$ of $2m$-tuples $(x_1, \ldots ,x_m: y_1, \ldots,
y_m)$ where $1 \leq x_1 < \cdots < x_m$ and where $y_1, \ldots ,y_m$
is order isomorphic to $p$.  In such a situation, we can say
$(p,\mathcal{S})$ occurs in a permutation $\pi = \pi_1 \ldots \pi_n$
if and only if there is a $2m$-tuple $(x_1, \ldots ,x_m: y_1,
\ldots, y_m) \in \mathcal{S}$ such that $\pi_{x_i} = y_i$ for $i=1,
\ldots, m$. While this is the most general type of pattern condition
that we can think of, in most cases this would be a very cumbersome
notation. Our definition of PDVP's was motivated by our attempts to
cover all the different types of pattern matching conditions that
have appeared in the literature that still allows for a relatively
compact notation.

\end{document}